\documentclass{cpaper}

\usepackage{amsfonts}
\usepackage[dvips]{graphics,color}
\newcommand{\Z}{\mathbb{Z}}
\newcommand{\Q}{\mathbb{Q}}

\newcommand{\R}{\mathbb{R}}

\input amssym.def
\newsymbol\square 1003
\newsymbol\lozenge 1006
\newsymbol\surjarrow 1310
\newsymbol\injarrow 131A
\newsymbol\ltimes 226E
\newsymbol\rtimes 226F
\newsymbol\void 203F
\newsymbol\leftrightsquigarrow 1321

\def\proclaim #1. #2\par #3\par {\medbreak
\noindent{\bf#1.\enspace}{\sl#2}\par\medbreak
\noindent{\bf Proof.} #3 \par}
\def\proclaimb #1. #2\par {\medbreak
\noindent{\bf#1.\enspace}{\sl#2}\par\medbreak}
 0
 2
 0
\newcommand{\arabiceqn}{\setcounter{equation}{0}%
\renewcommand{\theequation}{\mbox{\arabic{section}.\arabic{equation}}}}%
\newtheorem{thm}[equation]{Theorem}
\newtheorem{lemma}[equation]{Lemma}
\newtheorem{prop}[equation]{Proposition}
\newtheorem{cor}[equation]{Corollary}

\newtheorem{claim}[equation]{Claim}

\newtheorem{fact}[equation]{Fact}
\begin{document}
\title{Geometry for palindromic automorphism groups of free groups}
\author{Henry H. Glover and Craig A. Jensen}
\date{Winter 1999-2000}

\smallskip

\address{Department of Mathematics, The Ohio State University\\
Columbus, OH 43210, USA}
\email{glover@math.ohio-state.edu, jensen@math.ohio-state.edu}

\begin{abstract}
We examine the palindromic automorphism group $\Pi A(F_n)$ 
of a free group $F_n$, a group first defined by
Collins in \cite{[C2]} which is related to hyperelliptic 
involutions of mapping class groups, congruence subgroups of $SL_n(\Z)$,
and symmetric automorphism groups of free groups.  
Cohomological properties of the group are explored
by looking at a contractible space on which $\Pi A(F_n)$ acts 
properly with finite quotient. Our results answer
some conjectures of Collins and provide a few striking results
about the cohomology of $\Pi A(F_n)$, such as that its rational
cohomology is zero at the vcd.
\end{abstract}


\maketitle

\begin{center}
{\em Dedicated to Peter Hilton}: \\
{\sc doc, note, i dissent. a fast never 
prevents a fatness.  i diet on cod.} \\
\end{center}

\section{Introduction}
\arabiceqn

Let $Aut(F_n)$ be the automorphism group of a free group $F_n$ on $n$ 
generators $a_1, a_2, \ldots, a_n$. A reduced word
$x_1^{\epsilon_1}x_2^{\epsilon_2}\ldots x_n^{\epsilon_n}$ is
called a {\em palindrome} if it is equal to its reverse
$x_n^{\epsilon_n}x_{n-1}^{\epsilon_{n-1}}\ldots x_1^{\epsilon_1}$.
In \cite{[C2]} Collins defines the {\em palidromic automorphism group}
$\Pi A(F_n)$ as the subgroup of $Aut(F_n)$ consisting of all automorphisms
$\alpha$ for which $\alpha(a_i)$ is a palindrome for all $i$.  He
showed that the group was generated by three types of automorphisms:
\begin{itemize}
\item Maps $(a_i||a_j)$, $i \not = j$, which send $a_i \mapsto a_j a_i a_j$
and fix all other generators $a_k$.
\item Maps $\sigma_{a_i}$ which send $a_i \mapsto a_i^{-1}$ and fix all
other generators $a_k$.
\item Maps corresponding to elements of the symmetric group $\Sigma_n$
which permute the $a_1, \ldots, a_n$ among themselves.
\end{itemize}
The portion of $\Pi A(F_n)$ generated by just the $(a_i||a_j)$ is called
the {\em elementary palindromic automorphism group} of $F_n$ and
denoted $E\Pi A(F_n)$. Note that $\Pi A(F_n) = E\Pi A(F_n) \rtimes 
(\Z/2 \hbox{ } \wr \hbox{ }\Sigma_n).$  Collins showed that a set of defining 
relators for  $E\Pi A(F_n)$ is given by relations of the form
\begin{enumerate}
\item $(a_i||a_k)(a_j||a_k)=(a_j||a_k)(a_i||a_k)$
\item $(a_i||a_k)(a_j||a_l)=(a_j||a_l)(a_i||a_k)$
\item $(a_i||a_k)(a_j||a_k)(a_i||a_j)=
(a_i||a_j)(a_j||a_k)(a_i||a_k)^{-1}$
\end{enumerate}

He remarked how similar this was to the relations for the pure
symmetric automorphism group  $P\Sigma A(F_n)$ 
(see Gilbert's work in \cite{[G]}):
\begin{enumerate}
\item $(a_i|a_k)(a_j|a_k)=(a_j|a_k)(a_i|a_k)$
\item $(a_i|a_k)(a_j|a_l)=(a_j|a_l)(a_i|a_k)$
\item $(a_i|a_k)(a_j|a_k)(a_i|a_j)=
(a_i|a_j)(a_j|a_k)(a_i|a_k)$
\end{enumerate}
where $(a_i|a_j)$, $i \not = j$, sends $a_i \mapsto a_j^{-1} a_i a_j$
and fixes all other generators $a_k$.

On the basis of this, Collins conjectured that one could find the 
virtual cohomological dimension of $\Pi A(F_n)$ by employing the methods
of \cite{[C-V]}, as he did for $\Sigma A(F_n)$ in \cite{[C1]}.  
He also speculated that $E\Pi A(F_n)$ is torsion free, just as
 $P\Sigma A(F_n)$ is.  We are able to answer both of these questions
in this paper, as well as obtaining several interesting facts about the
cohomology of $\Pi A(F_n)$.

\begin{thm} \label{c34} Let $\Pi A(F_n)$ be the palindromic automorphism
group of the free group $F_n$ on $n$ letters and let $E\Pi A(F_n)$ be
the subgroup of elementary palindromic automorphisms. Then \newline
\noindent {\bf a) }
The virtual cohomological dimension of $\Pi A(F_n)$ is $n-1$. \newline
\noindent {\bf b) }
(i) For the prime 2, the Krull dimension of 
$\hat H^*(\Pi A(F_n); \Z_{(2)})$ is
$n$.  
For odd primes $p$, the Krull dimension of
$\hat H^*(\Pi A(F_n); \Z_{(p)})$ is $\left[\frac{n}{p}\right]$. \newline
(ii) 
In the range where the Krull dimension of $\hat H^*(\Pi A(F_n); \Z_{(p)})$
is 1, the period is $2(p-1)$. \newline
\noindent {\bf c) }
The group $E\Pi A(F_n)$ is torsion free. \newline
\noindent {\bf d) }
The cohomology group  $H^{n-1}(\Pi A(F_n); \Q) = 0$. \newline
\noindent {\bf e) }
If $p$ is an odd prime and $n = p, p+1, p+2$, then the
Farrell cohomology of the palindromic automorphism group is the 
same as that of the symmetric group on $p$ elements:
$$\hat H^*(\Pi A(F_n); \Z_{(p)}) \cong \hat H^*(\Sigma_p; \Z_{(p)}).$$
\end{thm}

For analogous results concerning $Aut(F_n)$, see \cite{[Ch]} and
\cite{[G-M]}.
See \cite{[F]} for the definition of the Farrell cohomology 
$\hat H^*(G; M)$ of a group $G$ of finite vcd with coefficients in
a $G$-module $M$ and also see \cite{[B]} for several useful properties of 
these cohomology groups.

The remainder of this paper is structured as follows.
In section 2, we discuss $\Pi A(F_n)$ and note how it relates to some other
groups, while in section 3 we introduce the space $L_{\sigma_n}$ which
$\Pi A(F_n)$ acts on and prove parts a) and d) of Theorem \ref{c34}.
Section 4 is concerned with a realization proposition which allows us
to establish parts b) (i) and c) of the main theorem.  Finally, section 5 
looks in more detail at the cohomology of $\Pi A(F_n)$ at odd primes $p$
and establishes parts b) (ii) and e) of the main theorem.


The authors would like to thank John Meier for an enlightening 
conversation about symmetric automorphisms of free groups, aiding
the presentation of this paper.

\section{Relationships with other groups}
\arabiceqn

Let $X_n$ be
the spine of auter space (see \cite{[C-V]},\cite{[H]},\cite{[H-V]}) 
and $Q_n = X_n/Aut(F_n)$.
Let $\sigma_n \in Aut(F_n)$ be
the automorphism which sends $a_i \mapsto a_i^{-1}$ for each $i$.  

Define the $\theta$-graph $\theta_m$ to be a graph with 2 vertices 
and $m+1$ edges, 
where each edge goes from one vertex to the other one.
Choose one of the two vertices of $\theta_1$ to be the basepoint $*$,
and define the rose $R_n$ to be the result of wedging together
$n$ copies of $\theta_1$ at the basepoint. 

The petals of the rose $R_n$
can be identified with the generators $a_i$ of $F_n$, so that 
$\pi_1(R_n,*) \cong F_n$. 
There is an action of $\langle\sigma_n\rangle = \Z/2$
on $R_n$ given by inverting each petal of the rose. This action
realizes the subgroup $\langle\sigma_n\rangle$ in the sense of $\cite{[Z]}$
(also cf. \cite{[C]}.)
An action of a group $G$ on a graph $\Gamma$ is {\em without inversions}
if $G$ does not send any edge $e$ to its inverse $\bar e$, and
an action is {\em reduced} if there are no $G$-invariant 
subforests in $\Gamma$.  The action of $\sigma_n$ on $R_n$ is both
without inversions and reduced.  From now on, when we refer to 
a group action on a graph, it is assumed that the edges of the
graph are subdivided as necessary to insure that the
group acts without inversions.

Note that the palindromic automorphism group $\Pi A(F_n)$ is just
$C_{Aut(F_n)}(\sigma_n)$.  This follows because an easy argument
shows that every element of $C_{Aut(F_n)}(\sigma_n)$ is palindromic,
and because the generators of  $\Pi A(F_n)$ are all in
$C_{Aut(F_n)}(\sigma_n)$.  For example, the generators $(a_i||a_j)$
are just products of $\sigma_n$-Nielsen transformations
(see \cite{[K]} where the $G$-Nielsen transformation
$\langle e, f \rangle \Gamma$ of a $G$-graph
$\Gamma$ has the same vertex and edge set as $\Gamma$ but where the
terminal point of an edge $eg$, $g \in G$, in the new graph is
the initial point of the edge $fg$ in the original graph; this induces
a map $\langle e, f \rangle$ from the fundamental groupoid of the
first graph to that of the second where $eg$ is sent to $(ef)g$ and all
other edges are sent to themselves.)
That is, if the petal $a_i$ of the rose 
consists of the edges $\bar e_i f_i$, then   
$(a_i||a_j)$ is the composition $\langle e_i, \bar f_j\rangle 
\circ \langle e_i, e_j \rangle$.

As a note for the curious, it follows that $\Pi A(F_n)$ and
$\Sigma A(F_n)$ are distinct groups for $n \geq 2$, 
since a direct argument shows that $\Sigma A(F_n)$ has no 
element of order 2 in its center.  
In addition, $E\Pi A(F_n)$ and $P\Sigma A(F_n)$ are also obviously
distinct (for $n \geq 3$; for $n=1,2$ they are the same group with
the same presentation,) since the former abelianizes to an
elementary abelian $2$-group of rank $n(n-1)$ while the latter 
abelianizes to a free abelian group of rank $n(n-1)$.

In addition to its formal palindromic properties,
the group $\Pi A(F_n)$ arises naturally from looking at hyperelliptic
subgroups of mapping class groups (cf. Gries \cite{[GR]} for corresponding
homological properties.)
We have a commutative diagram
\begin{equation} \label{c24}
\matrix{\Gamma_g^{2,pure} & \rightarrow & \Gamma_g^1 
& \rightarrow & \Gamma_g \cr
\downarrow & & \downarrow & & \downarrow \cr
Aut(F_{2g}) & \rightarrow & Out(F_{2g}) & \rightarrow & GL_{2g}(\Z) \cr}
\end{equation}
where $\Gamma_g$ is the mapping class group of an orientable surface of
genus $g$, $\Gamma_g^1$ is the mapping class group of an orientable 
surface of genus $g$ with 1 puncture, and $\Gamma_g^{2,pure}$ is the 
mapping class group of an orientable surface of genus $g$ with two punctures,
where each puncture is fixed pointwise. The map from $\Gamma_g^{2,pure}$
to $Aut(F_{2g})$ is obtained first by taking an intersection basis 
$a_1, b_1, \ldots, a_g, b_g$ for the fundamental group of the surface $S$
with two punctures.  One of the punctures should serve as the basepoint for
fundamental group considerations.  The other is treated as an actual puncture,
so that the fundamental group of the surface minus this point is a free group
$F_{2g}$ on $2g$ generators.  The map from $\Gamma_g^{2,pure}$
to $Aut(F_{2g})$ is now obtained by sending an element of
$\Gamma_g^{2,pure}$ to the automorphism of $F_{2g}$ that it induces.
The map from  $\Gamma_g^1$ to $Out(F_{2g})$ is obtained similarly.

Let $\psi \in \Gamma_g^{2,pure}$ be a hyperelliptic involution 
(see, for example, \cite{[F-K]}.)  Then $\psi$ has $2g+2$ fixed points, 
two of which are of course the punctures on the surface $S$.  Choose loops
$a_1, b_1, \ldots, a_g, b_g$, based at one of the punctures, which form
an intersection basis for the surface $S$ (say the ones described 
in \cite{[F-K]} in the section on hyperelliptic Riemann surfaces.)  
By going along the top row of diagram \ref{c24} and then projecting downward,
we see that the image of $\psi$ in $GL_{2g}(\Z)$ is $-I$. Let $\bar \psi$ be
the image of $\psi$ in $Aut(F_{2g})$.  Our goal is to show that
$\bar \psi$ is conjugate to $\sigma_{2g}$ in $Aut(F_{2g})$.

\begin{lemma} \label{c25}
If $\phi \in Aut(F_n)$ and  the image of $\phi$ in $GL_n(\Z)$ is $-I$,
then the image of any conjugate $\alpha^{-1}\phi\alpha$ of $\phi$,
$\alpha \in Aut(F_n)$, is also $-I$.
\end{lemma}
\begin{proof} This follows directly since $-I$ is in the center of $GL_n(\Z)$. \end{proof} 

\begin{lemma} \label{c26}
If $\phi \in Aut(F_n)$ is an involution whose image in $GL_n(\Z)$ is $-I$,
then $\phi$ can be realized on a marked graph whose underlying graph is
the rose.
\end{lemma}

\begin{proof} Realize $\phi$ on a reduced marked graph whose underlying graph is
$\Lambda$.  First, we show that if 
$e$ is an edge of $\Lambda$ which is not fixed by $\phi$,
then we can assume that one endpoint of $e$ is the basepoint $*$.
Let $f = \phi(e)$.  Choose a shortest path $\gamma$ from $e$ to $*$.
Since $stab(e)= \langle 1 \rangle$, $stab(e) \subseteq stab(h)$ for every $h$
in the path $\gamma$.  So we can apply a sequence of 
Nielsen transformations (see \cite{[K]})
and slide $e$ along $\gamma$ to $*$.  Note that since $\Lambda$ is reduced,
$\{e,f\}$ now forms either a rose $R_2$ based at $*$, or a $\theta_1$.
Proceeding in this manner, we can slide all of the edges of $\Lambda$ not
fixed by $\phi$ to the basepoint.

By way of contradiction, suppose that two of these $\theta_1$-graphs
(that have been moved so that one vertex of each $\theta_1$ 
is the basepoint),
say $\{e_1,f_1\}$ and $\{e_2,f_2\}$, share another common vertex in
addition to the basepoint $*$.  That is, suppose that there is a vertex
$v \not = *$ and each of $e_1,f_1,e_2,f_2$ go from $v$ to $*$.
Say $g: R_n \to \Lambda$ is the marked graph, and recall that
$\pi_1(R_n) = \langle a_1, \ldots, a_n \rangle = F_n$.  By replacing
$\phi$ by a conjugate if necessary (see Lemma \ref{c25}) we can assume
that $g$ sends the petal $a_1$ of $R_n$ to $e_1^{-1}f_1$,
the petal $a_2$ to $e_2^{-1}f_2$, and the petal
$a_3$ of $R_n$ to $e_1^{-1}f_2$. So in $\pi_1(\Lambda)$ we have
$$\phi \cdot a_1=f_1^{-1}e_1=g(a_1^{-1}),$$
$$\phi \cdot a_2=f_2^{-1}e_2=g(a_2^{-1}),$$
and
$$\phi \cdot a_3 = f_1^{-1}e_2= g(a_1^{-1}a_3a_2^{-1}).$$
Hence the first column of $im(\phi)$ in $GL_n(\Z)$ is $(-1,0,\ldots,0)$,
the second column is $(0,-1,0,\ldots,0)$,
and the third column is $(-1,-1,1,0,\ldots,0)$.  This contradicts the fact
that $im(\phi)=-I$.

If the result of sliding to the basepoint $*$ all edges of $\Lambda$ not fixed
by $\phi$ yields a rose $R_n$, then we are done.
Otherwise,
suppose by way of contradiction  
that there exist edges $e$,$f$, and $h$ of
$\Lambda$ such that \newline
\begin{itemize}
\item Both $e$ and $f$ go from some vertex $v \not = *$ to $*$.
\item $\phi(e)=f$.
\item $h$ goes from $v$ to $v$.
\item $\phi(h)=h$.
\end{itemize}

As before, say $g: R_n \to \Lambda$ is the marked graph. By replacing
$\phi$ by a conjugate if necessary (see Lemma \ref{c25}) we can assume
that $g$ sends the petal $a_1$ of $R_n$ to $e^{-1}f$ and the petal
$a_2$ of $R_n$ to $e^{-1}hf$.  So in $\pi_1(\Lambda)$ we have
$$\phi \cdot a_1=f^{-1}e=g(a_1^{-1})$$ and
$$\phi \cdot a_2 = f^{-1}he = (e^{-1}f)^{-1}(e^{-1}hf)(e^{-1}f)^{-1}
= g(a_1^{-1}a_2a_1^{-1}).$$
Hence the first column of $im(\phi)$ in $GL_n(\Z)$ is $(-1,0,\ldots,0)$
and the second column is $(-2,1,0,\ldots,0)$.  This contradicts the fact
that $im(\phi)=-I$. \end{proof}

\begin{prop} \label{c27}
If $\phi \in Aut(F_n)$ is an involution whose image in $GL_n(\Z)$ is $-I$,
then $\phi$ is conjugate in $Aut(F_n)$ to $\sigma_n$.
\end{prop}

\begin{proof} Realize $\phi$ on a marked graph $g: R_n \to R_n$.  By replacing
$\phi$ by a conjugate if necessary, we can assume $g(a_i)=a_i$ for all
$i$.  The involution $\phi$ of the graph $R_n$ must send the petal
$a_1$ to some $a_j^{\pm 1}$, since it is a graph automorphism.  But since
$im(\phi) = -I \in GL_n(\Z)$, we see that $a_j^{\pm 1}$ must be $a_1^{-1}$.
Similarly, we can see that the graph automorphism $\phi$ sends, for each $i$,
the petal $a_i$ to the petal $a_i^{-1}$.  This means our current
$\phi$ is equal to $\sigma_n$ (and thus our original $\phi$, before 
we replaced it by a conjugate, was conjugate to $\sigma_n$.) \end{proof}

The following corollary is immediate:

\begin{cor} \label{c28}
The image $\bar \psi$ of $\psi$ in $Aut(F_{2g})$ is conjugate to 
$\sigma_{2g}$.  Hence the hyperelliptic subgroup 
$C_{\Gamma_g^{2,pure}}(\psi)$ is conjugate to a subgroup of
$C_{Aut(F_{2g})}(\sigma_{2g})$.
\end{cor}

As remarked in \cite{[C2]}, the image of $\Pi A(F_n)$ in 
$GL_n(\Z)$ is the subgroup of $GL_n(\Z)$ consisting of invertible
matrices where each column has exactly one odd entry (and the rest are 
even.)  The subgroup is
the semidirect product $\tilde \Gamma_2(\Z)  \rtimes \Sigma_n$ 
where $\tilde \Gamma_2(\Z)$ is 
the 2-congruence subgroup
defined by the short exact sequence
$$ \langle 1 \rangle \to \tilde \Gamma_2(\Z) \injarrow GL_n(\Z) 
\surjarrow GL_n(\Z/2) \to \langle 1 \rangle$$
and $\Sigma_n$ is standard inclusion of the symmetric group
$$\Sigma_n \subset GL_n(\Z).$$

\section{A space for $\Pi A(F_n)$ to act on}
\arabiceqn

We define a certain contractible space $L_{\sigma_n}$, 
related to auter space, which
$\Pi A(F_n)$ acts on with finite stabilizers and finite quotient.  This
allows us to obtain some cohomological results.

A graph $\Gamma$ is a {\em $\theta_1$-tree of rank $n$} if there
exists a pointed tree $T$ such that $\Gamma$ is obtained by 
``doubling'' every edge of $T$ into a $\theta_1$-graph.  That is, the
vertex set of $\Gamma$ is the same as the vertex set of $T$ and for every
edge $e$ of $T$ going from $v$ to $w$, $\Gamma$ has two edges
$e_1$ and $e_2$, both of which go from $v$ to $w$.  There is a
natural $\Z/2$-action on such a graph $\Gamma$, which is given by
switching the two edges in each $\theta_1$-graph.  Note that
the orbit space of $\Gamma$ under this action is just the tree $T$.

\begin{claim} \label{c1} The reduced graphs $\Gamma$ which realize
the subgroup $\langle\sigma_n\rangle$ of $Aut(F_n)$ are exactly the
$\theta_1$-trees of rank $n$, where $\sigma_n$ acts on the trees
via their natural $\Z/2$-action.
\end{claim}

\begin{proof}  We have already mentioned that the rose $R_n$ realizes
$\sigma_n$.
From Theorem 2 of \cite{[K]}, the other reduced graphs 
$\Gamma$ which also realize $\sigma$ are those that are Nielsen 
equivalent to $R_n$ (up to an equivariant isomorphism.)

If $e$ is an edge in one of the copies of $\theta_1$ in $R_n$
and $f$ is an edge in a different $\theta_1$ in $R_n$, and both
$e$ and $f$ point toward the basepoint, then Nielsen transformation
$\langle e,f\rangle$ has the result of pulling the $\theta_1$-graph
$\{e, \sigma_n e\}$ through the $\theta_1$-graph
$\{f, \sigma_n f\}$, so that now $e$ terminates at
the initial vertex of $f$, rather than at the basepoint $*$.
In other words,
the result of applying one Nielsen transformation to $R_n$ is that
of sliding one of the petals of $R_n$ up though another petal.
 
A basic induction argument now yields that the result of applying 
a series of Nielsen transformations fo $R_n$ will be some 
$\theta_1$-tree $\Gamma$. \end{proof}

Recall from \cite{[K-V]} that an edge $e$ of a $G$-graph $\Gamma$ is
{\em inessential} if it is in every maximal $G$-invariant subforest of
$\Gamma$.  A $G$-graph $\Gamma$ is {\em inessential} if it has at
least one inessential edge and is {\em essential} if it is not 
inessential.  Let $X_n^G$
be the fixed point subspace of $X_n$ corresponding to some finite
subgroup $G$ of $Aut(F_n)$.  From \cite{[J2]} (cf. part III of \cite{[JD]}
and \cite{[K-V]}), both
the centralizer $C_{Aut(F_n)}(G)$ and the 
normalizer $N_{Aut(F_n)}(G)$ act on the contractible space $X_n^G$ with
finite stabilizers and finite quotient.  Moreover, the space
$X_n^G$ $G$-equivariantly deformation retracts to the space $L_G$,
where $L_G$ is constructed from $X_n^G$ by considering only
essential marked graphs.  Hence $L_G$ is a good space to study if one 
wishes to calculate the cohomology of $C_{Aut(F_n)}(G)$
or $N_{Aut(F_n)}(G)$. 

Further recall that a $G$-graph $\hat \Gamma$ is a $G$-equivariant
{\em blowup} of a $G$-graph $\Gamma$ if some $G$-invariant subforest $F$ of
$\hat \Gamma$ can be collapsed away to yield $\Gamma$. 
Let $\phi: R_n \to \Gamma$ be some reduced marked graph realizing 
$\langle\sigma_n\rangle$.  From Claim \ref{c1}, $\Gamma$ is a $\theta_1$-tree.
Blow up $\Gamma$ $\sigma_n$-equivariantly to some maximal essential 
blowup $\hat \Gamma$.

\begin{claim} \label{c2} The fixed points/cells of the 
action of $\sigma_n$ on $\hat \Gamma$ are exactly the valence 2 vertices 
of $\hat \Gamma$.
\end{claim}

\begin{proof} Note that vertices in $\hat \Gamma$ 
have valence 2 or 3 and that 
$\hat \Gamma$ is obtained from $\Gamma$ by blowing up an oriented 
ideal forest (see \cite{[J2]}, \cite{[K-V]}, \cite{[JD]}.)  Briefly, ideal
edges (oriented ideal forests) correspond to subsets of edges 
(chains of subsets of edges) which are pulled away from existing vertices
in order to create new graphs which collapse down to the original graph.)

No edge is in $Fix_{\sigma_n}(\hat \Gamma)$ because no edge is in
$Fix_{\sigma_n}(\Gamma)$ and so blowing up ideal edges will not create
any new edges that are fixed under the action of $\sigma_n$ (see \cite{[K-V]}
page 229.)

If a valence 3 vertex is in $Fix_{\sigma_n}(\hat \Gamma)$, then at least
one edge of $\hat \Gamma$ must be fixed by $\sigma_n$, which is a
contradiction.

All of the valence 2 vertices of $\Gamma$ are in $Fix_{\sigma_n}(\Gamma)$.  
New valence 2 vertices created as ideal edges are blown up correspond 
to either:
\begin{itemize}
\item Old vertices of $\Gamma$ that used to be valence higher than 2 but 
have since had edges stripped (pulled away) from them.  These are in
$Fix_{\sigma_n}(\hat \Gamma)$ because all vertices of $\Gamma$ are in
$Fix_{\sigma_n}(\Gamma)$.
\item New valence 2 vertices inserted to insure that $\sigma_n$ acts on
$\hat \Gamma$ without inversions.  These are also clearly in
$Fix_{\sigma_n}(\hat \Gamma)$. \end{itemize}  \end{proof}

Note that $* \in Fix_{\sigma_n}(\Gamma)$.  Cut $\hat \Gamma$ along each of
its valence 2 vertices, yielding a graph $\hat \Gamma_{cut}$ with the
same number of valence 3 vertices and edges as $\hat \Gamma$ had, no valence
2 vertices, and twice as many valence 1 vertices as $\hat \Gamma$ had valence
2 vertices.

\begin{claim} \label{c3} 
$$\hat \Gamma_{cut} = \hat \Gamma_{1} \amalg \hat \Gamma_{1},$$
the disjoint union of two trees $\hat \Gamma_{1}$ and $\hat \Gamma_{2}$,
where $\hat \Gamma_{2} = \sigma_n \hat \Gamma_{1}$.
\end{claim}

\begin{proof} There is a covering map $p: \hat \Gamma_{cut} \to \hat \Gamma/\sigma_n$
obtained by mapping to the orbit space under the $\sigma_n$-action.  
The forest collapse that sends $\hat \Gamma$ to $\Gamma$ is 
$\sigma_n$-equivariant, so it descends to a forest collapse of
$\hat \Gamma/\sigma_n$ to $\Gamma/\sigma_n$.  (That the
quotient of the forest upstairs in $\hat \Gamma$ is also a forest in
$\hat \Gamma/\sigma_n$ can be seen by an easy Euler characteristic
argument.) But $\Gamma$ is a $\theta_1$-tree with a known $\sigma_n$-action
on it, and $\Gamma/\sigma_n$ is a tree (in fact, the underlying tree of
the $\theta_1$-tree.)  Hence $\hat \Gamma/\sigma_n$ is a tree.  Since
$p: \hat \Gamma_{cut} \to \hat \Gamma/\sigma_n$ is a covering map with
fiber two points, $\hat \Gamma_{cut}$ is as described. \end{proof}

Let $T$ be a pointed tree with $2n-1$ edges, all vertices valence 
either 1 or 3, where $*$ is one of the valence 1 vertices.  (Then $T$ has
$n+1$ valence 1 vertices and $n-1$ valence 3 vertices.)  Let $T_1$ and
$T_2$ be two isomorphic copies of $T$, and let $f: T_1 \to T_2$
be an isomorphism.  Define
$$\Gamma_T = \frac{T_1 \amalg T_2}{f(v) \sim v, 
\hbox{ for all valence 1 vertices $v$ of $T_1$.}}$$
Define a $\sigma_n$-action on $\Gamma_T$ by
$$\sigma_n x = \left\{ \matrix{f(x), \hfill &x \in T_1 \hfill \cr
f^{-1}(x), \hfill &x \in T_2 \hfill \cr} \right.$$

\begin{prop} \label{c4} There is a bijective correspondence between 
trees $T$ as above (with $2n-1$ edges, etc) and maximal, essential 
blowups $\hat \Gamma$ of reduced $\sigma_n$-graphs.  The bijection is
given by $T \mapsto \Gamma_T$.
\end{prop}

\begin{proof} From claims \ref{c2} and \ref{c3}, all blowups $\hat \Gamma$ have the
required form.  Finally, any $\Gamma_T$ can easily be reduced to a 
$\theta_1$-tree by collapsing edges, meaning that it is the blowup of
such a graph. \end{proof}

All maximal simplices in $L_{\sigma_n}$ have the same dimension,
from \cite{[K-V]}.  Maximal simplices in 
$L_{\sigma_n}/\Pi A(F_n)$ are constructed by taking 
chains of forest collapses from maximal blowups $\Gamma_T$.  Alternatively,
we can define a {\em subforest} of $T$ to be a collection $S$ of
edges of $T$ such that there is no path in $S$ from one valence 1 vertex
to another.  (If there were such a path, then $S \cup \sigma_n S$ would
be a cycle in $\Gamma_T$.)  In this way, we can think of maximal simplices
as coming from chains of subforests of various trees $T$.

\medskip

\begin{proof}[Proof of part a) of Theorem \ref{c34}:] Since
$\Pi A(F_n)$ acts on the
contractible
space $L_{\sigma_n}$ with finite stabilizers and finite quotient,
the vcd of $\Pi A(F_n)$ is at most the
dimension of a maximal cell from $L_{\sigma_n}/\Pi A(F_n)$.
Such a cell comes from a chain of forest collapses of a tree $T$
with $2n-1$ edges, $n+1$ valence 1 vertices, and $n-1$ valence 3
vertices.  Hence we can collapse at most $n-1$ of the valence 3
vertices into other vertices while doing forest collapses, resulting
in maximal simplices of dimension $n-1$.

To show that the vcd of  $\Pi A(F_n)$ is at least $n-1$,
we note that the subgroup generated by $(a_i || a_n)$ for 
$i \in \{1, 2, \ldots, n-1\}$  is isomorphic to
$\Z^{n-1} = \Z \times \ldots \times \Z$. \end{proof}
 
\begin{lemma} \label{c6} Let $\tilde \Gamma$ be the underlying
graph of a particular marked graph in $L_{\sigma_n}$.  Hence
$\tilde \Gamma$ comes equipped with a $\sigma_n$-action.  Let $C$
be a simple closed curve in $\tilde \Gamma$.  Then $Fix_{\sigma_n}(C)$
contains exactly two points, and $\sigma_n(C)$ is the curve $-C$, or
$C$ with the opposite of its original orientation.
\end{lemma}

\begin{proof} $\tilde \Gamma$ can be blown up (not necessarily uniquely) to
some maximal essential graph $\Gamma = \Gamma_T$.  $\Gamma$ is the
union of two isomorphic copies $T_1$ and $T_2$ of $T$, where $T_1$ and
$T_2$ are attached along their corresponding valence 1 vertices.

As we collapse from $\Gamma$ to $\tilde \Gamma$, the trees $T_1$ and
$T_2$ collapse to trees $\tilde T_1$ and $\tilde T_2$.  However, 
the attaching points for $\tilde T_1$ and $\tilde T_2$ are no longer 
necessarily just the valence 1 vertices, and could be other vertices as
well.

Let $\alpha_1$ be a taut path in $\tilde T_1$ from one attaching 
point $v_1$ to some other attaching point $v_2$, where 
furthermore there are no attaching points in the interior of $\alpha_1$.
Let $\alpha_2 = \sigma_n(\alpha_1)$ be the corresponding path in 
$\tilde T_2$.  Then $\alpha_1 \bar \alpha_2$ is a simple closed 
curve and $\sigma_n(\alpha_1 \bar \alpha_2$) = $\alpha_2 \bar \alpha_1$,
or the original curve oriented in the other direction.  Hence
the curve $\alpha_1 \bar \alpha_2$ satisfies the conclusions of the 
lemma.  Our goal is to show that any simple closed curve $C$ takes 
this form.

Let $C_1$ be the portion of $C$ that is in $\tilde T_1$ and let
$C_2$ be the portion of $C$ that is in $\tilde T_2$. Since $C$ is a cycle
and yet both $\tilde T_1$ and $\tilde T_2$ are trees, both $C_1$ and $C_2$
are nonempty.  In fact, there must be a path 
$\alpha_1$ in $C_1$ from one attaching point $v_1$ to some other 
attaching point $v_2$, where there are no other attaching points in the
interior of the path.  Since $C$ is a simple closed curve, $\alpha_1$
must be a taut path.  Let $\alpha_2 = C - \alpha_1$, some other taut
path in $\tilde \Gamma$ from $v_1$ to $v_2$.

Let $p : \tilde \Gamma \to \tilde T_1$ be the map given by taking the
quotient space under the action of $\sigma_n$.  Note that $p(\alpha_2)$
is a path from $v_1$ to $v_2$. Since $\tilde T_1$ is a tree and
$\alpha_1$ is the unique taut path in $\tilde T_1$ from $v_1$ to $v_2$,
this gives us that 
$\hbox{edges}(\alpha_1) \subseteq \hbox{ edges}(p(\alpha_2))$.  Hence if
$e$ is an edge in $\alpha_1$, then $p^{-1}(e)$ is two edges,
$e \in C_1$ and $\sigma_n(e) \in C_2$. It follows that all of the 
oriented edges of the simple closed curve 
$\alpha_1 \sigma_n(\bar \alpha_1)$ are in the simple closed curve $C$.
Hence $C = \alpha_1 \sigma_n(\bar \alpha_1).$ \end{proof}

\begin{lemma} \label{c7} Let $\tilde \Gamma$ be the underlying
graph of a particular marked graph in $L_{\sigma_n}$.  Hence
$\tilde \Gamma$ comes equipped with a $\sigma_n$-action. 
Choose an edge $e$ in $\tilde \Gamma$.
Among all simple closed
curves $D$ which pass through $e$, choose one curve $C$ for which the 
distance from the curve to the basepoint $*$ is minimal.  
Then $Fix_{\sigma_n}(C)$ is two points $v_1$ and $v_2$.
One of these two points in the closest point in $C$ to the basepoint
$*$ and one of them is the farthest point in $C$ to the basepoint $*$.
\end{lemma}

\begin{proof} Using the notation of the proof of Lemma \ref{c6}, $C$ is the result
of following some path $\alpha_1$ in $\tilde T_1$ and then
$\bar \alpha_2 = \sigma_n (\bar \alpha_1)$ in $\tilde T_2$, where
$\alpha_1$ goes from the attaching point $v_1$ to the attaching point $v_2$.

Since $\alpha_1$ is a path in a pointed tree, there is a unique vertex
$w$ in $\alpha_1$ which is closest to $*$.
By way of contradiction, suppose that $w \not \in \{v_1, v_2\}$. 
Then $w$ is not an attaching point.
Let $\beta_1$ be the unique taut path in $\tilde T_1$ from $*$ to $w$.
Let $\gamma_1$ be the unique subpath of $\beta_1$ which contains $w$ and
exactly one attaching point $y$. Now let $\delta_1$ be the path in
$\tilde T_1$ which  starts at $y$, follows $\gamma_1$ along to $w$, and
then either follows $C$ from $w$ to $v_2$ or $-C$ from $w$ to $v_1$
(where we choose whichever possibility insures that $\pm e \in \delta_1$.)
Then the simple closed curve $\delta_1 \sigma_n(\bar \delta_1)$ is closer
to $*$ than $C$ is, which is a contradiction.
Hence $w$ is $v_1$ or $v_2$, and the lemma follows. \end{proof}

\begin{prop} \label{c8}  Let $\tilde \Gamma$ be a graph which 
occurs as an underlying graph  of marked graphs in $L_{\sigma_n}$. 
Then there is only one possible $\sigma_n$-action on $\tilde \Gamma$.
\end{prop}

\begin{proof} We see that our task is to show that a 
unique $\sigma_n$-action is determined by the properties about
simple closed curves listed in Lemmas \ref{c6} and \ref{c7}.

Define an action $\eta$ on $\tilde \Gamma$ as follows.
Let $e$ be an oriented edge of $\tilde \Gamma$.  Among all simple closed
curves $D$ which pass through $e$, choose one path $C$ for which the 
distance from the curve to the basepoint $*$ is minimal. Let $v_1$ be
a point on $C$ which is closest to the basepoint.  Let $n$ be the 
edge-path distance in $C$ from $v_1$ to $e$.  Then there is an orientation
$\epsilon \in \{-1,1\}$ such that if you traverse $\epsilon C$ starting
at $v_1$ and go $n$ edges, you get to $e$.  Define $\eta(e)$ to be the
result of traversing $-\epsilon C$ starting at $v_1$ and then going $n$
more edges.

By Lemma \ref{c6} and \ref{c7}, the action $\eta$ is well defined and
if any $\sigma_n$ acts on $\tilde \Gamma$, then the $\sigma_n$-action and the
$\eta$-action coincide. \end{proof}

Denote by $Q_{\sigma_n}$ the quotient space
$L_{\sigma_n}/\Pi A(F_n)$. 

\begin{cor} \label{c9} If two marked graphs in $L_{\sigma_n}$ have the
same underlying graph, then they correspond to the same vertex in 
$Q_{\sigma_n}$.  That is, the moduli space $Q_{\sigma_n}$ can be formed
by looking only at the poset structure of the underlying graphs of
marked graphs in $L_{\sigma_n}$.
\end{cor}

\begin{proof} From Proposition \ref{c8}, any underlying graph of a
marked graph in $L_{\sigma_n}$ has only 
one possible $\sigma_n$-action.  But from Corollary 10.4 of \cite{[K-V]},
$\Pi A(F_n)$ acts transitively on the set of marked 
$\sigma_n$-graphs based on the same $\sigma_n$-graph.  The result
follows. \end{proof}

The simplices of $L_{\sigma_n}$ group themselves into cubes, as described
in \S 3 of \cite{[V]}.  In \cite{[H-V]}, Hatcher and Vogtmann
show that the quotients in $Q_n$ of cubes in $X_n$ have the rational 
homology of balls.  They use this to create a cubical chain complex
which has the same rational homology as $Q_n$.  Our goal here is to 
establish a similar result for the cubes of maximal dimension in
$Q_{\sigma_n}$, where this time we want the quotients of cubes to have 
the $\Z_{(p)}$-cohomology of balls, where $p$ is any odd prime.

Following \cite{[H-V]}, we consider a maximal cube in $L_{\sigma_n}$.  It is
given by considering a maximal essential marked graph 
$\phi : R_n \to \Gamma_T$ and considering some maximal subforest $S$ of $T$.
Recall that by a subforest of $T$ we mean a subset $S$ of the edges of 
$T$ where there is no path in $S$ from one valence 1 vertex of $T$
(or equivalently, one terminal edge of $T$) to another.
From part a) of Theorem \ref{c34}, 
$S$ has $n-1$ edges in it.  The cube corresponding
to the pair $(T,S)$ can thus be thought of as imbedded in $\R^{n-1}$, where
each coordinate vector is an edge of the cube, the graph obtained by
collapsing each edge of $S \cup \sigma_n(S)$ is at the origin,
and $\Gamma_T$ is at $(1,1,\ldots,1)$.  Let $Aut(T,S)$ be 
the group of all (pointed) automorphisms of the tree $T$ which take
$S$ to $S$. The group 
$stab_{\Pi A(F_n)}(T,S) = \langle\sigma_n\rangle \times Aut(T,S)$ acts 
linearly on the cube by permuting the coordinates of $\R^{n-1}$, and
fixes the diagonal from $(0,0,\ldots,0)$ to $(1,1,\ldots,1)$.  (The
involution $\sigma_n$ acts trivially on the cube, of course, since all
of these cubes are coming from $L_{\sigma_n} \subset X_n^{\sigma_n}$.)
Hence, just as in \cite{[V]}, the quotient of the cube in 
$Q_{\sigma_n}$ is a cone with base $S^{n-2}/Aut(T,S)$, where $S^{n-2}$ is
the boundary of the cube.

\begin{lemma} \label{c10} The finite group $Aut(T)$ 
(and hence its subgroup $Aut(T,S)$) is all $2$-torsion.
\end{lemma}

\begin{proof} Let $\xi \in Aut(T)$.  Now $\xi$ must take the basepoint to the
basepoint, and so it must take the unique edge attached to the basepoint
to itself.  For each $n$, let $E_n$ be the edges in $T$ which are at most
distance $n$ from $*$.  So $E_0$ is just one edge, and $\xi$ fixes it 
as already mentioned.  Since all nonterminal vertices of $T$ have
valence 3, an inductive argument yields that
$\xi^{2^n}$ fixes $E_n$ pointwise. \end{proof}

The following is an analog of Proposition 3.1 of \cite{[V]}:

\begin{prop} \label{c11} $S^{n-2}/Aut(T,S)$ has the $\Z_{(p)}$-cohomology
of an $(n-2)$-sphere or a ball.  The latter possibility happens when there
is an element of $Aut(T,S)$ which induces an odd permutation of the
edges of $S$.
\end{prop}

\begin{proof} The finite group $Aut(T,S)$, which is all 2-torsion, acts cellularly 
on $S^{n-2}$, where the stabilizer of a cell fixes it pointwise.
We use the spectral sequence for equivariant cohomology 
(cf. \cite{[B]} VII \S 7):  
\begin{equation} \label{c12}
E_1^{r,s} = \prod_{[\delta] \in \Delta_n^r}
H^s(stab(\delta); \Z_{(p)}) \Rightarrow H_{Aut(T,S)}^{r+s}(S^{n-2}; \Z_{(p)})
\end{equation}
where $[\delta]$ ranges over the set $\Delta^r$ of orbits
of $r$-simplices $\delta$ in $S^{n-2}$.
Since $Aut(T,S)$ is all $2$-torsion and finite, so are all of the
$stab(\delta)$.  Hence if $s > 0$, $H^s(stab(\delta); \Z_{(p)})=0$.
So the above spectral sequence converges to \newline
$H^r(S^{n-2}/Aut(T,S); \Z_{(p)})$.

But another filtration yields a spectral sequence with
\begin{equation} \label{c13}
E_2^{r,s} = H^r(Aut(T,S); H^s(S^{n-2}; \Z_{(p)}))
\Rightarrow H_{Aut(T,S)}^{r+s}(S^{n-2}; \Z_{(p)})
\end{equation}
It follows that $E_2^{r,s} = 0$ unless $(r,s)$ is $(0,0)$ or $(0,n-2)$.
Hence $E_2^{0,0} = \Z_{(p)}$ and 
$E_2^{0,n-2} = H^{n-2}(S^{n-2}; \Z_{(p)})^{Aut(T,S)}$.  The latter
group of invariants is $\Z_{(p)}$ if the action of $Aut(T,S)$ on
$S^{n-2}$ preserves orientation and $0$ otherwise.  The last assertion
in the proposition follows from Corollary 3.2 of \cite{[V]}. \end{proof}

\begin{thm} \label{c14} The top dimensional cohomology group of 
$Q_{\sigma_n}$ vanishes.  That is, $H^{n-1}(Q_{\sigma_n};\Z_{(p)}) = 0$.
\end{thm}

\begin{proof} We show that the quotient of  every maximal cube 
$(T,S)$ has a free face, so that the interior of the quotient of the cube 
can be collapsed away.  If we can do this, then $Q_{\sigma_n}$ will
have the same $\Z_{(p)}$-cohomology as an $(n-2)$-dimensional complex,
and we will be done.

In the degenerate case where there is an element of $Aut(T,S)$ which 
induces an odd permutation of the edges of $S$, then the quotient of the
cube $(T,S)$ is not itself a cube.  
In this case, the diagonal from $(0,\ldots,0)$
to $(1,\ldots,1)$ is exposed in the quotient, and any $(n-2)$-dimensional
simplex in the quotient which lies next to the diagonal is a free face.

In the nondegenerate case, the quotient of the cube $(T,S)$ is itself a
cube, although its boundary might be self indentified in various ways.
Since the subforest $S$ of $T$ is maximal, $S$ must contain at least one
terminal edge $e$.  That is, one of the two vertices of $e$ is a valence
1 vertex or attaching point.  Let $\tilde \Gamma$ be the graph obtained
from $\Gamma = \Gamma_T$ by collapsing the subforest $\{e,\sigma_n(e)\}$.
The graph $\tilde \Gamma$ has a maximal subforest corresponding to 
collapsing the edges $e$ and $\sigma_n(e)$ from the forest
$S \cup \sigma_n(S)$ of $\Gamma$.  Hence we see that collapsing $e$ gives
us a face, which we will denote by $(T/e,S/e)$, of the cube $(T,S)$.
It can be shown that this 
face corresponds to a (nondegenerate, cubical) face of the
quotient of the cube $(T,S)$ because

\begin{claim} \label{c15} There is an natural injection of 
$Aut(\tilde \Gamma) = \langle\sigma_n\rangle \times Aut(T/e)$ into
$Aut(\Gamma) = \langle\sigma_n\rangle \times Aut(T)$.  Define the lift
$\hat \phi$ of an automorphism $\phi \in Aut(T/e)$
by sending an edge $f$ to $\phi(f)$ if
$f \not = e$ and letting $\hat \phi(e) = e$. 
\end{claim}

\begin{proof} Denote by $v$ the valence 1 vertex of $e \in T$ 
(the attaching point) and let $w$ be the other vertex of $e$.  
In $T/e$, $w=v$. We must show that $\phi$ sends $w$ to $w$.  
This follows automatically, however, as $w=v=\sigma_n(w)=\sigma_n(v)$ 
is the only valence 4 vertex of $\tilde \Gamma$ and so any 
automorphism of the graph must fix it.  Let $f$ and $g$ be the two other 
edges in $T$ which share the vertex $w$.  Now if $v=*$ then $\phi$ could
possibly exchange $f$ and $g$, but this is fine as the lift
$\hat \phi$ also can.  If $v \not = *$, then one of $f$ or $g$ must
be closer to the basepoint, and so $\phi$ must fix both $f$ and $g$.
Regardless, $\hat \phi$ can be defined as in the statement of the claim.
\end{proof}

\noindent Warning:  Note that if $e$ is not a terminal edge, the
above claim is false.  Collapsing an interior edge somtimes allows
you to construct automorphisms with 3-torsion, which obviously
cannot be lifted to $Aut(T)$.

No automorphism $\phi$ of $(T/e,S/e)$ can 
induce an odd permutation of the edges in $S/e$,
else the lift $\hat \phi$ of $\phi$ to $T$ would
induce an odd permutation of the edge of $S$.
Since $Aut(T)$ is all 2-torsion, it follows from the above claim that
$Aut(T/e)$ is also all 2-torsion.   Hence the same 
spectral sequence argument used in Proposition \ref{c11} yields that
the quotient of the cube corresponding to $(T/e,S/e)$
actually is a $\Z_{(p)}$-cohomology cube.

It remains to be shown that the cubical face corresponding to $(T/e,S/e)$
is free.  First, if another subforest $S'$ with an edge $e'$ of $T$ 
gives a cube with a face isomorphic to $(S/e,S/e)$, then $e'$ must also 
be a terminal edge of $S'$.  Hence the isomorphism  
$(T/e,S/e) \to (T/e',S'/e')$ maps the vertex that $e$ collapsed into to the
the vertex that $e'$ collapsed into, and so we can lift the 
isomorphism to one from $(T,S) \to (T,S')$.
 
Second, we must show that blowing up the vertex $w$ in $\tilde \Gamma$
only yields graphs isomorphic to $\Gamma$.
This follows by considering the ways that the vertex $w$ in
$\tilde \Gamma$ can be blown up. Say that the edges $f$, $g$,
$\sigma_n(f)$, and $\sigma_n(g)$ are the ones incident to $w$.  
If the ideal edge orbit  $\sigma_n\{f,g\}$ is blown up, we get
back $\Gamma$ exactly, and if $\sigma_n\{f,\sigma_n(g)\}$ is blown up, 
we get a graph isomorphic to $\Gamma$.  As these are the only ways to
blow up the graph $\sigma_n$-equivariantly into another essential
graph, we are done. \end{proof}

\begin{cor} \label{c16}  
$H^{n-1}(Q_{\sigma_n};\Q) = H^{n-1}(\Pi A(F_n); \Q) = 0$.
\end{cor}

\begin{proof} That $H^{n-1}(Q_{\sigma_n};\Q) = 0$ follows immediately from
Theorem \ref{c14}.  Recall that $\Pi A(F_n)$ acts with
finite stabilizers and finite quotient $Q_{\sigma_n}$ on the 
contractible space $L_{\sigma_n}$.  Since the stabilizers are finite, 
their rational cohomology vanishes, and the standard equivariant
spectral sequence yields that 
$H^{*}(Q_{\sigma_n};\Q) = H^{*}(\Pi A(F_n); \Q)$. \end{proof}

Note that part d) of Theorem \ref{c34} follows from the above 
Corollary.

As a final remark for this section, we show that $L_{\sigma_n}$ is
an $\underline{E} \Pi A(F_n)$ (cf. \cite{[K-M]}); 
that is, for finite subgroups $G$
of $\Pi A(F_n)$, the fixed point subcomplex $L_{\sigma_n}^G$ is 
contractible.  This follows directly from the corresponding property
of $Aut(F_n)$.  The following proposition is unneccesary in
the specific case of $L_{\sigma_n}$, since (proof omitted) $L_{\sigma_n}$ 
actually equals $X_n^{\sigma_n}$.  This does not normally happen
(for example, the spaces $L_{P_n \times \sigma_n}$ mentioned later in
Fact \ref{c21} are not equal to the corresponding fixed point 
space of $X_n$), however, and thus it seems worth noting the more
general fact.

\begin{prop} \label{ptwo1} Let $S$ be a finite subgroup of $Aut(F_n)$ and
let $\mathcal{S}$ be either $C_{Aut(F_n)}(S)$ or $N_{Aut(F_n)}(S)$.  Let
$L_S$ be the retract, defined by Krstic and Vogtmann and consisting of
essential marked graphs, of the fixed point subcomplex $X_n^S$ of the spine
of auter space $X_n$.  Then $L_S$ is an $\underline{E} \mathcal{S}$ space.
\end{prop}

\begin{proof}[Sketch of Proof:]
Let $H$ be a finite subgroup of $\mathcal{S}$ and let $G$ be the
(finite, because $HSH^{-1}=S$) subgroup generated by $H$ and $S$.  
Then $X_n^G = (X_n^S)^H = (X_n^H)^S$, and $X_n^G$
is contractible from \cite{[J2]}.  It remains to be shown that 
$X_n^G = (X_n^S)^H$ deformation retracts to $(L_S)^H$.

Given a marked graph $\Gamma$ representing a vertex of $(X_n^S)^H$, 
we must show (see Proposition 3.3 of \cite{[K-V]}) that for
every edge $e$ in $\Gamma$ and every $h \in H$, $e$ is $S$-inessential
if and only if $he$ is $S$-inessential.  This follows automatically
from Corollary 4.5 of \cite{[K-V]}, which characterizes 
essential edges by looking 
at the stabilizers (in $S$) of paths in $\Gamma$.  Since $HSH^{-1}=H$, 
the stabilizers in $h$-translates of such paths are still in $S$ 
and are isomorphic (conjugate by $h$) to those of the original path.
\end{proof}

\section{A realization proposition}
\arabiceqn

Let $\hat A$ be a finite subgroup of $\Pi A(F_n)$ and let $A$ be the
(finite) subgroup generated by $\hat A$ and $\sigma_n$.
By Zimmerman's \cite{[Z]} realization theorem, we can realize 
$A$ by an action on an $A$-reduced graph $\Gamma$.  From the 
proposition below, $\Gamma$ is also $\langle\sigma_n\rangle$-reduced; 
that is, $\Gamma$ is a $\theta_1$-tree.  

Note that the corresponding statement is not true
in $Out(F_n)$ (have $\Z/p \times \langle\sigma_{p-1}\rangle$ 
act on a $\theta$-graph
$\theta_{p-1})$ and certainly would not be true in $Aut(F_n)$ if the
$\sigma_n$-action were replaced by some other $\Z/2$-action. 

\begin{prop} \label{c18} 
Let $A \subseteq \Pi A(F_n)$ be a finite subgroup of the
palindromic automorphism group with $\sigma_n \in A$.  
Realize $A$ by an action on an
$A$-reduced marked graph $\phi: R_n \to \Gamma$.  Then $\phi: R_n \to \Gamma$
is also a $\langle\sigma_n\rangle$-reduced marked graph.
\end{prop}

\begin{proof} As before, let $F_n = \langle a_1, \ldots, a_n\rangle$ and identify the petals
of the rose $R_n$ with the generators $a_i$.  Note that $\Gamma$ has
no separating edges, else it would not be $A$-reduced. In this proof,
when we refer to concepts such as the number of times an edge $e$ of 
$\Gamma$ occurs in some $\phi(a_i)$, we mean that we should take the
unique taut path in $\Gamma$, starting and ending at $*$, 
which is homotopic to the path $\phi(a_i)$ in $\Gamma$, and then 
count the number of times $e$ occurs in this taut path.
By way of contradiction, suppose $\Gamma$ is not $\langle\sigma_n\rangle$-reduced.
Let $e_1 \in \Gamma$ be an edge of minimal 
distance to the basepoint $*$ such that
$\{e_1,\sigma_n e_1\}$ is a forest.

\noindent {\bf CASE 1: } $e_1 = \sigma_n e_1$.  Since $e_1$ is not a
separating edge of $\Gamma$, we can choose a nontrivial cycle $\mu$, starting
and ending at $*$, which has just one occurence of $e_1$ and none of 
$e_1^{-1}$.  If for all $i=1, \ldots, n$, the cycles $\phi(a_i)$ have an
even number of occurences of $e_1^{\pm 1}$, then we could not write $\mu$
as a product of them and their inverses.  So some $\phi(a_j)$ has an 
odd number of occurences of $e_i^{\pm 1}$ in it.  Say that the exponent
sum of $e_1$ in $\phi(a_j)$ is $k$, $k$ odd.  Then the exponent sum of
$e_1$ in $\sigma_n \phi(a_j)$ is still $k$, but the exponent sum of
$\phi(a_j^{-1})$ is $-k$.  This contradicts the fact that 
$\sigma_n a_j = a_j^{-1}.$

\noindent {\bf CASE 2: } $e_1 \not = \sigma_n e_1$. Let $\alpha$ be a
shortest length path  from $*$ to $e_1$.  Say without loss of generality
that $e_1$ is the oriented edge from $v$ to $w$ and that $\alpha$ goes
from $*$ to $v$.  Let $f_1 = \sigma_n e_1$.  Then $\sigma_n \alpha$ is a
shortest length path from $*$ to $f_1$.  Now $v=\sigma_n v$, else we
could write $\alpha = \beta b$ and get $\{b, \sigma_n b\}$ as a 
$\sigma_n$-invariant forest closer to $*$ than $\{e_1,f_1\}$ is.
(If $|\alpha|=0$, then $v=*$ and so $\sigma_n(v)=v$ necessarily.)
So we have both $\alpha$ and $\sigma_n \alpha$ are paths from $*$ to $v$ 
Moreover,  $w \not = \sigma_n w$ 
(else $\{e_1,\sigma_n e_1\}$ is not a forest.)
Now $A e_1 = A \{e_1, f_1\}$ is not a forest, since $\Gamma$ is $A$-reduced.
Hence we can choose some simple closed curve $\mu$ in 
$A e_1 \subseteq \Gamma$ that contains $e_1$.  There must exist some
$a e_1^{\pm 1} \in \mu$, $a e_1^{\pm 1} \not \in \{e_1, e_1^{-1}\}$,
such that $a w = w$.  Why?  Otherwise we could deformation retract
$\mu$ to the set of vertices $\{ \hat a v: \hat a e_1^{\pm 1} \in \mu \}$,
which contradicts the fact that $\mu$ is a simple closed curve.
Now $a e_1 \not = f_1$, as $\sigma_n w \not = w$.  Hence
$\alpha e_1 (a e_1)^{-1} (a \alpha)^{-1}$ is a nontrivial cycle starting
and ending at $*$ which contains exactly one occurence of $e_1$ and
none of $f_1$.  So there must be a $\phi(a_j^{\epsilon})$,
$\epsilon \in \{-1, 1\}$,  which contains an
odd number of occurences of $e_1^{\pm 1}$ and an even number of occurences
of $f_1^{\pm 1}$.  (If we had some even/odd $\phi(a_j)$, then we could act by
$\sigma_n$ to get odd/even, and this would be a $\phi(a_j^{-1})$. 
Otherwise, all $\phi(a_i)$ are all even/even or
odd/odd, and so combine together just to get more even/even or odd/odd
loops.) This is a contradiction, however, because $\phi(a_j^{-\epsilon})$
still has an odd number of occurences of  $e_1^{\pm 1}$ while
$\phi(\sigma_n a_j^\epsilon)$ has an even number of
occurences of $e_1^{\pm 1}$. \end{proof}

\medskip

\begin{proof}[Proof of part b) (i) of Theorem \ref{c34}:] 
From the action of $(\Z/2)^n$ on the rose $R_n$
such that the $i$th generator inverts the $i$th petal and leaves
all others fixed,
we know 
that the Krull dimension at the prime 2 is at least $n$.  
Similarly, there is an action of
$(\Z/p)^{\left[\frac{n}{p}\right]}$ on $R_n$ where the first
$\Z/p$ 
rotates the first $p$ petals, the second $\Z/p$ rotates the
next $p$ petals, etc.  Hence  the Krull dimension 
at the prime $p$ is at least
$\left[\frac{n}{p}\right].$

Let $A$ be a maximal rank elementary abelian subgroup of $\Pi A(F_n)$.
From Proposition \ref{c18}, we can realize $A$ by an action of $A$ on a 
$\sigma_n$-graph $\Gamma$ which is both $A$-reduced and $\sigma_n$-reduced.
That is, we have an action of $A$ on a pointed $\Theta_1$-tree $\Gamma$.
Since elements of $A$ must preserve basepoints,
the action of $A$ on the tree $\Gamma/\sigma_n$ does not invert edges.
Hence we have inclusions
$$A \injarrow (\Z/2)^n \rtimes Aut_*(\Gamma/\sigma_n) 
\injarrow (\Z/2)^n \rtimes \Sigma_n = \Z/2 \wr \Sigma_n.$$
The result (for 2 or odd primes $p$) now follows from standard facts
about $\Sigma_n$ (cf. Theorem 1.3 in Chapter VI of \cite{[A-M]}.)
\end{proof}

\medskip

\begin{proof}[Proof of part c) of Theorem \ref{c34}:]
We sketch the proof, which uses standard methods.
Suppose that some $A = \Z/p$ lies in $E\Pi A(F_n)$. From
Proposition \ref{c18},  we can realize $A$ by an action of $A$ on a 
$\sigma_n$-graph $\Gamma$ which is both $A$-reduced and 
$\sigma_n$-reduced.  Let $\phi: R_n \to \Gamma$ be the corresponding
marked graph. Let $T=\Gamma/\sigma_n$, a pointed tree with an
$A$-action on it.  First, suppose that the $A$-action on $T$ is
nontrivial.  (This will always be the case if $p$ is odd.)  Then 
there are two
edges $e_1$ and $e_2$, both oriented so that their terminal vertices are
closer to the basepoint than their initial vertices, 
of $\Gamma$ such that a generator of $A$ rotates the edge
$[e_1]$ into the edge $[e_2]$ in $T$.   Some generator $a_{i_1}$ of
$F_n$ must be such that $\phi(a_{i_1})$ contains an odd number of occurences
of $e_1$ in it.  Choose $a_{i_2}$ similarly.  Then $\phi(a_{i_j})$ is
a palindromic word in the edges of $\Gamma$ with either 
$\bar e_{i_j} \sigma_n(e_{i_j})$ or $\sigma_n(\bar e_{i_j}) e_{i_j}$
in the middle of the palindrome.  The generator of $A$ (thought of
as an element of $E\Pi A(F_n)$) must send $a_{i_1}$ to a palindrome with
either $a_{i_2}$ or $a_{i_2}^{-1}$ in the center of it.  This 
contradicts the fact that all elements of $E\Pi A(F_n)$ send generators
$a_i$ to palidromes with $a_i$ in the center of them.

The only remaining case is where $p=2$ and $A$ acts trivally on $T$.
So $A$ is a subgroup of the group $(\Z/2)^n$ of graph automorphisms
of $\Gamma$ which act by inverting the $\Theta_1$'s in the $\Theta_1$-tree
$\Gamma$.  Hence the generator of $A$ corresponds to an element of 
$\Pi A(F_n)$ which, for at least one $i$, sends $a_i$ to a 
palindrome with $a_i^{-1}$ in its center.  As none of these automorphisms
are in $E\Pi A(F_n)$, we again have a contradiction. \end{proof}


\section{Cohomology of $\Pi A(F_n)$ at odd primes $p$}
\arabiceqn

Let $p$ be an odd prime (as will always be the case
from now on in this paper.)  
We wish to calculate the Farrell cohomology of $\Pi A(F_n)$ using 
Ken Brown's \cite{[B]} normalizer spectral sequence, which states that
\begin{equation} \label{c17}
E_1^{r,s} = \prod_{(P_0 \subset \cdots \subset P_r) \in |\mathcal{B}|_r}
\hat H^s( \bigcap_{i=0}^r N_{G}(P_i); \Z_{(p)})
\Rightarrow \hat H^{r+s}(G; \Z_{(p)})
\end{equation}
where $G$ is a group with finite virtual cohomological dimension,
$\mathcal{A}$ is the poset of nontrivial elementary abelian
$p$-subgroups of $G$,
$\mathcal{B}$ is the poset of conjugacy classes of
nontrivial elementary
abelian $p$-subgroups of $G$, and $|\mathcal{B}|_r$ is the
set of $r$-simplices in the realization $|\mathcal{B}|$.

A first step toward performing such a calculation is calculating
$|\mathcal{B}|$.  In other words, we wish to calculate conjugacy classes
of elementary abelian subgroups $P \subset \Pi A(F_n)$. 
By Proposition \ref{c18},  we can realize such finite groups 
$P$ by reduced actions on $\theta_1$-trees.

If $n \geq p$, define a particular subgroup $P_n \cong \Z/p$ of 
$\Pi A(F_n)$ by letting $P_n$ act on the rose $R_n$ by 
rotating its first $p$ leaves and leaving the last $n-p$ leaves fixed.
That is, $P_n$ corresponds to automorphisms which rotate the first
$p$ generators $a_1, \ldots, a_p$ and leave the remaining generators fixed.

\begin{cor} \label{c19} If $p \leq n \leq 2p-1$, then 
$$\hat H^*(\Pi A(F_n); \Z_{(p)}) \cong
\hat H^*(N_{Aut(F_n)}(P_n \times \langle\sigma_n\rangle); \Z_{(p)}).$$
\end{cor}

\begin{proof} We show that $P_n$ is the only conjugacy class of nontrivial
elementary abelian $p$-subgroups that is in $\Pi A(F_n)$.
By Proposition \ref{c18}, we see that an arbitrary nontrivial elementary
abelian $p$-subgroup $A$ comes from some action on a $\theta_1$-tree
with $p$-symmetry.  Since $p \leq n \leq 2p-1$, the only possibility is
that $A$ acts on a $\theta_1$-tree $\Gamma$ by rotating $p$ of the 
$\theta_1$-leaves and leaving the other $n-p$ $\theta_1$-edges in the
tree fixed.  But it is clear that a product of 
$(P_n \times \langle\sigma_n\rangle)$-Nielsen transformations takes the rose $R_n$
to the graph $\Gamma$, and hence we see that $A$ and $P_n$ are 
conjugate to each other in $\Pi A(F_n)$.

By the normalizer spectral sequence \ref{c17}, this yields that
$$\hat H^*(\Pi A(F_n); \Z_{(p)}) \cong
\hat H^*(N_{\Pi A(F_n)}(P_n); \Z_{(p)}).$$
But since $p$ is an odd prime, it is easy to see that
$$N_{\Pi A(F_n)}(P_n) = 
N_{Aut(F_n)}(P_n \times \langle\sigma_n\rangle).$$ \end{proof}

\begin{prop} \label{c20}
$$N_{Aut(F_n)}(P_n \times \langle\sigma_n\rangle) \cong
N_{\Sigma_p}(P_n) \times 
(F_m \rtimes (\langle\sigma_p\rangle \times \Pi A(F_m)))$$  
where $m = n-p$, $\Pi A(F_m)$ acts on the $F_m$ in the semidirect product
in the natural way, and $\sigma_p$ acts on $F_m$ as
$\sigma_m$ does.
\end{prop}

\begin{proof} The $N_{\Sigma_p}(P_n)$ in the above decomposition comes from
automorphisms of $F_n$ which permute the first $p$ generators and leave
the remaining $m$ fixed.  The $F_m$ being acted upon in the semidirect
product structure above has $i$th generator 
$(a_1||a_{p+i})(a_2||a_{p+i}) \ldots (a_p||a_{p+i})$.  The $\sigma_p$ is
the involution which inverts the first $p$ generators of $F_n$ and leaves
the remaining $m$ fixed.  Finally, the $\Pi A(F_m)$ comes from
automorphisms which fix the first $p$ generators of $F_n$ and act on the
last $m$ generators by
identifying the subgroup $\langle a_{p+1}, a_{p+2}, \ldots, a_n \rangle$
with $F_m$.

Consider the action of $P_n \times \langle\sigma_n\rangle$ on the rose $R_n$.
$P_n$ rotates the first $p$ petals. Label the first $p$ petals of the
rose as $a_1, \ldots, a_p$ as before, but label the last $m$ petals
as $b_1, \ldots, b_m$.

Since $|P_n|=p$ is an odd prime, $N_{Aut(F_n)}(P_n \times \langle\sigma_n\rangle) \subseteq 
N_{Aut(F_n)}(P_n)$ and in Lemma 5.1 of \cite{[J]}, we calculated
$$N_{Aut(F_n)}(P_n) \cong 
N_{\Sigma_p}(P_n) \times 
((F_m \times F_m) \rtimes (\langle\sigma_p\rangle \times Aut(F_m))),$$
where the first $F_m$ in $F_m \times F_m$ is the free group
on the $P_n$-Nielsen transformations $\langle a_1,b_i^{-1}\rangle$ for
$i \in \{1, \ldots, m\}$ and the latter $F_m$ is the free
group on the $P_n$-Nielsen transformations 
$\langle a_1^{-1}, b_i^{-1}\rangle$, $i \in \{1, \ldots, m\}$. 
Note that $\langle\sigma_p\rangle$ acts on $F_m \times F_m$ via
$\sigma_p \langle a_1,b_i^{-1}\rangle \sigma_p = \langle a_1^{-1}, b_i^{-1}\rangle$ and
$\sigma_p \langle a_1^{-1},b_i^{-1}\rangle \sigma_p = \langle a_1,b_i^{-1}\rangle$.
In other words, if $(b,c) \in F_m \times F_m$ then
$\sigma_p (b,c) \sigma_p = (c,b)$.

Let $G$ be the 
subgroup 
$$N_{\Sigma_p}(P_n) \times 
(F_m \rtimes (\langle\sigma_p\rangle \times C_{Aut(F_m)}(\sigma_m)))$$
of $N_{Aut(F_n)}(P_n)$, where $F_m$ is the free group on the
generators $\langle a_1,b_i\rangle \circ \langle a_1^{-1},b_i^{-1}\rangle$ 
for 
$i \in \{1, \dots, m\}$,
and $C_{Aut(F_m)}(\sigma_m)$ is included in $Aut(F_m)$ 
in the obvious way.  It follows directly that 
$G \subseteq N_{Aut(F_n)}(P_n \times \langle\sigma_n\rangle).$  To
prove the proposition, we must show that they are equal.

Take an arbitrary
$$x \in N_{Aut(F_n)}(P_n \times \langle\sigma_n\rangle) \subseteq 
N_{\Sigma_p}(P_n) \times ((F_m \times F_m) \rtimes 
(\langle\sigma_p\rangle \times Aut(F_m))).$$
Say $x = abcde$, where $a \in N_{\Sigma_p}(P_n)$,
$(b,c) \in F_m \times F_m$, $d \in \langle\sigma_p\rangle$,
and $e \in Aut(F_m)$. Since $a, d \in N_{Aut(F_n)}(P_n \times \langle\sigma_n\rangle)$,
$a^{-1}xd^{-1}=bce \in  N_{Aut(F_n)}(P_n \times \langle\sigma_n\rangle)$ also.
So $bce \in  \Pi A(F_n)$ and
$(bce) \sigma_n (bce)^{-1} = \sigma_n.$  This means that the
map $(bce) \sigma_n (bce)^{-1}$ sends $a_i$ to $a_i^{-1}$ for
$i \in \{1,\ldots,p\}$ and $b_i$ to $b_i^{-1}$ for
$i \in \{1,\ldots,m\}$.  Now both $\sigma_n$ and $e$ restrict to
maps in $Aut(\langle b_1, \ldots, b_m\rangle)$ and moreover $b$ and $c$ both
restrict to the identity map in  $Aut(\langle b_1, \ldots, b_m\rangle)$.
Hence for $i \in \{1,\ldots,m\}$, we have
$$b_i^{-1} = (bce) \sigma_n (bce)^{-1} (b_i) = e \sigma_n e^{-1} (b_i),$$ 
and we see that $e \sigma_n e^{-1}$ restricts to $\sigma_m$ in
$Aut(F_m)$. As $e \in Aut(F_m)$, this means
$e \in C_{Aut(F_m)}(\sigma_m)$.  Hence $e \in \Pi A(F_n)$
also.  Since $bce \in \Pi A(F_n)$, this gives
$bc \in  \Pi A(F_n)$.  In other words, we have
$$(b,c) \in (F_m \times F_m) \subseteq N_{\Sigma_p}(P_n) \times ((F_m \times F_m) \rtimes (\langle\sigma_p\rangle \times Aut(F_m)))$$
and
$$(b,c) \in \Pi A(F_n).$$
It follows that
$$\matrix{\hfill (b,c) &=& \sigma_n (b,c) \sigma_n \hfill \cr
\hfill &=& \sigma_m \sigma_p (b,c) \sigma_p \sigma_m \hfill \cr
\hfill &=& \sigma_m (c,b) \sigma_m \hfill \cr
\hfill &=& (\sigma_m(c),\sigma_m(b)). \hfill \cr}$$
So $b=\sigma_m(c)$ and $c=\sigma_m(b)$.  In summary, we have shown that
an arbitrary element $x=abcde \in 
N_{Aut(F_n)}(P_n \times \langle\sigma_n\rangle)$ has $c = \sigma_m(b)$
and $e \in C_{Aut(F_m)}(\sigma_m)$.
Thus $x \in G$, as desired. \end{proof}

The group $N_{Aut(F_n)}(P_n \times \langle\sigma_n\rangle)$ acts on the
contractible space $L_{P_n \times \langle\sigma_n\rangle}$ with finite stabilizers
and finite quotient $Q_{P_n \times \langle\sigma_n\rangle}= 
L_{P_n \times \langle\sigma_n\rangle}/N_{Aut(F_n)}(P_n \times \langle\sigma_n\rangle)$.

Define a {\em $p$-admissible tree} $T$ to be a triple $(T,\circ,A)$ where
$T$ is a pointed tree, $\circ$ is a vertex of $T$ (which may be the basepoint
$*$), $A$ is a subset of the vertices of $T$ called the
{\em set of attaching points}, $* \in A$, and all valence 1
vertices of $T$ are in $A$.  For a $p$-admissible tree $T$, define the
{\em corresponding graph} $\Gamma_T$ as follows:  Take two isomorphic
copies $T_1$ and $T_2$ of the tree $T$, and let $f: T_1 \to T_2$ be an
isomorphism.  Then let $\Gamma_T^{pre}$ be the graph
$$\Gamma_T^{pre} = \frac{T_1 \amalg T_2}{f(v) \sim v, 
\hbox{ for all attaching points $v$ in $A$.}}$$
Let $\theta_{p-1}$ be a $\theta$-graph with $p$ edges and
two vertices $v_1$ and $v_2$.  Let $\circ_1$ be the $\circ$-vertex
in $T_1$ and let $\circ_2 = f(\circ_1)$ be the $\circ$-vertex in $T_2$. 
Finally, let
$$\Gamma_T = \frac{\Gamma_T^{pre} \amalg \theta_{p-1}}
{\circ_1 \sim v_1, \circ_2 \sim v_2}$$
If $\pi_1(\Gamma_T) \cong F_n$, then say $T$ is a $p$-admissible tree
of {\em rank} $n$.

If $T$ is a $p$-admissible tree of rank $n$,
define a $\langle\sigma_n\rangle$-action on the edges of $\Gamma_T$ by
$$\sigma_n x = \left\{ \matrix{f(x), \hfill &x \in T_1 \hfill \cr
f^{-1}(x), \hfill &x \in T_2 \hfill \cr
x^{-1}, \hfill &x \in \theta_{p-1} \hfill \cr} \right.$$
Since this action inverts the edges of the $\theta$-graph in
$\Gamma_T$, we then need to subdivide these edges so that the
group acts without inversions.  Next, define a $P_n$-action on $\Gamma_T$
by having $P_n$ fix $\Gamma_T^{pre}$ and 
rotate the edges of $\theta_{p-1}$ cyclically.  In this way,
$\Gamma_T$ is a $(P_n \times \langle\sigma_n\rangle)$-graph.

A $p$-admissible tree is $T$ {\em reduced} if the corresponding
$(P_n \times \langle\sigma_n\rangle)$-graph $\Gamma_T$ is reduced; that is, if all
vertices of $T$ are attaching points.  Similarly, a $p$-admissible tree $T$
is a {\em maximal} if the attaching points of $T$ are 
exactly its valence 1 vertices, the valence $2$ vertices of $T$ consist
of just the point $\circ$, and $T$ has no vertices with valence 4 or more.
As before, a {\em subforest} of $T$ is a collection of edges 
$S$ of $T$ such that there is no path in $S$ from 
one attaching point to another.  Lastly, isomorphisms of $p$-admissible
trees must be graph isomorphisms which 
take $*$ to $*$, $\circ$ to $\circ$, and $A$ to $A$.

The following facts about $(P_n \times \langle\sigma_n\rangle)$-graphs are all proven 
in similar ways to the analogous facts about $\sigma_n$-graphs.

\begin{fact} \label{c21}
\begin{enumerate}
\item There is a bijective correspondence between reduced $p$- admissible
trees of rank $n$ and the underlying graphs of $(P_n \times \langle\sigma_n\rangle)$-
reduced marked graphs, given by $T \to \Gamma_T$.
\item There is a bijective correspondence between maximal $p$-admissible
trees of rank $n$ and the underlying graphs of maximal essential marked
$(P_n \times \langle\sigma_n\rangle)$-graphs, given by $T \to \Gamma_T$.
\item The virtual cohomological dimension of 
$N_{Aut(F_n)}(P_n \times \langle\sigma_n\rangle)$ is $m=n-p$.
\item Let $\Gamma$ be a graph which occurs as the underlying graph of
a marked graph in $L_{P_n \times \langle\sigma_n\rangle}$.  Then there is only one 
possible $\sigma_n$-action on $\Gamma$.
\item  If two marked graphs in $L_{P_n \times \langle\sigma_n\rangle}$ have 
underlying graphs which correspond to the same $p$-admissible tree, 
then they correspond to the same vertex in 
$Q_{P_n \times \langle\sigma_n\rangle}$.  That is, we can form 
the moduli space $Q_{P_n \times \langle\sigma_n\rangle}$
by looking only at the poset structure of the $p$-admissible trees
corresponding to marked graphs in $L_{P_n \times \langle\sigma_n\rangle}$.
\item The top dimensional cohomology class of 
$Q_{P_n \times \langle\sigma_n\rangle}$, with 
coefficients in $\Z_{(p)}$, vanishes.  
That is, $H^{n-p}(Q_{P_n \times \langle\sigma_n\rangle};\Z_{(p)}) = 0$.
\item $H^{n-p}(Q_{P_n \times \langle\sigma_n\rangle};\Q) = 
H^{n-p}(N_{Aut(F_n)}(P_n \times \langle\sigma_n\rangle); \Q) = 0$.
\end{enumerate}
\end{fact}

Note that (4) and (5) above are a little bit different from their
analogs Proposition \ref{c8} and Corollary \ref{c9}.  Basically,
the underlying graphs $\Gamma$ always have just one possible 
$\sigma_n$-action, as before, but it is conceivable (for example, if
the graph contains two or more copies of $\theta_{p-1}$ inside it and
we must decide which one $P_n$ rotates) that it might have several
possible $P_n$-actions.  That is why we talk about $p$-admissible 
trees instead in (5), since the vertex $\circ$ in the tree determines
where the $p$ edges that $P_n$ rotates are located.

Fact \ref{c21} allows us to show 

\begin{prop} \label{c22} If $p \leq n \leq 2p-1$, then 
$$\hat H^{t}(N_{Aut(F_n)}(P_n \times \langle\sigma_n\rangle); \Z_{(p)}) 
\cong \left\{\matrix{
\Z/p \hfill &t \equiv 0  \hbox{ } (\hbox{mod } 2(p-1)) \hfill \cr
H^r(Q_{P_n \times \langle\sigma_n\rangle}; \Z/p) \hfill 
&t \equiv r  \hbox{ } (\hbox{mod } 2(p-1)), \hfill \cr
 &1 \leq r \leq n-p-1 \hfill \cr
0 \hfill 
&t \equiv r  \hbox{ } (\hbox{mod } 2(p-1)), \hfill \cr
 &n-p \leq r \leq 2p-3 \hfill \cr
} \right.$$
\end{prop}

\begin{proof} We use the equivariant cohomology spectral sequence for \newline
$N_{Aut(F_n)}(P_n \times \langle\sigma_n\rangle)$ acting on the contractible
space $L_{P_n \times \langle\sigma_n\rangle}$ with finite stabilizers and
finite quotient $Q_{P_n \times \langle\sigma_n\rangle}$. 
The equivariant cohomology spectral sequence
for this action is
$$\matrix{
E_1^{r,s} = \prod_{[\delta] \in \Delta_{n}^r}
\hat H^s(stab_{N_{Aut(F_n)}(P_n \times \langle\sigma_n\rangle)}(\delta); 
\Z_{(p)}) \hfill \cr
\hfill \Rightarrow 
\hat H^{r+s}(N_{Aut(F_n)}(P_n \times \langle\sigma_n\rangle)); \Z_{(p)}) \cr}$$
where $[\delta]$ ranges over the set $\Delta_{n}^r$ of orbits
of $r$-simplices $\delta$ in $L_{P_n \times \langle\sigma_n\rangle}$.

From the decomposition
$$N_{Aut(F_n)}(P_n \times \langle\sigma_n\rangle) \cong
N_{\Sigma_p}(P_n) \times 
(F_m \rtimes (\langle\sigma_p\rangle \times C_{Aut(F_m)}(\sigma_m)))$$  
we see that $(F_m \rtimes (\langle\sigma_p\rangle \times C_{Aut(F_m)}(\sigma_m)))$
has $p$-rank 0.  Since $N_{\Sigma_p}(P_n)$ acts 
trivially on marked graphs in $L_{P_n \times \langle\sigma_n\rangle}$ by 
permuting the edges of the $\theta$-graph attached at $\circ$,
it follows that for every simplex $\delta$ we have
$$\hat H^*(stab_{N_{Aut(F_n)}(P_n \times \langle\sigma_n\rangle)}(\delta); \Z_{(p)})
 \cong \hat H^*(N_{\Sigma_p}(P_n); \Z_{(p)}) \cong 
\hat H^*(\Sigma_p; \Z_{(p)}).$$

The $E_1^{r,s}$-page of the spectral sequence is $0$ in the rows
where $s \not = k \cdot 2(p-1)$ and a copy
of the cellular cochain complex with $\Z/p$-coefficients of the
$(n-p)$-dimensional complex $Q_{P_n \times \langle\sigma_n\rangle}$ in rows 
$k\cdot2(p-1)$.  It follows that
the $E_2$-page has the form:
$$E_2^{r,s} = \left\{\matrix{
\Z/p \hfill &r=0 \hbox{ and } s=k\cdot2(p-1) \hfill \cr
H^r(Q_{P_n \times \langle\sigma_n\rangle}; \Z/p) \hfill &1 \leq r \leq n-p
\hbox{ and } s=k\cdot2(p-1) \hfill \cr
0 \hfill &\hbox{otherwise } \hfill \cr} \right.$$
Hence we see that the spectral sequence converges at the $E_2$-page.

That $H^{n-p}(Q_{P_n \times \langle\sigma_n\rangle}; \Z/p) = 0$ follows from 
part 6 of Fact \ref{c21} and universal coefficients. \end{proof}

Note that the above proposition immediately proves
part b) (ii) of Theorem \ref{c34}.

\input{hyper1.pic}
\begin{figure}[here]
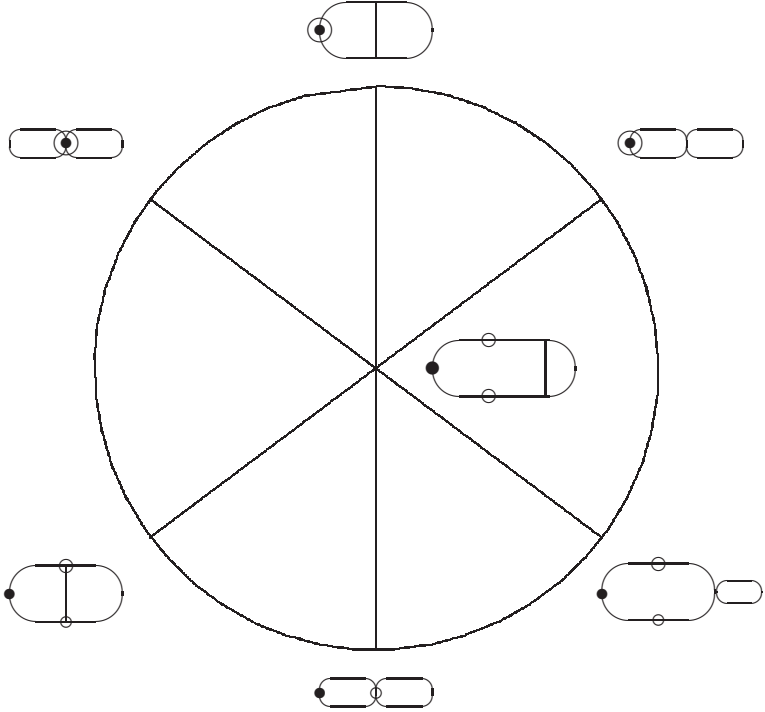

\caption{\label{hyper1} Simplices from the first maximal graph}
\end{figure}

By examining the space $Q_{P_n \times \langle\sigma_n\rangle}$ 
in low dimensions where 
$m \in \{0,1,2\}$ and showing that it is contractible, we have
the following corollary, which will give us part e) of 
Theorem \ref{c34}:

\begin{cor} \label{c23}
If $m = n-p \in \{0,1,2\}$, then
$$\hat H^*(\Pi A(F_n); \Z_{(p)}) \cong
\hat H^*(N_{Aut(F_n)}(P_n \times \langle\sigma_n\rangle); \Z_{(p)}) \cong
\hat H^*(\Sigma_p; \Z_{(p)}) .$$
\end{cor}

\begin{proof} {\bf CASE 1:} m=0.  Then $Q_{P_n \times \langle\sigma_n\rangle}$ 
is a point. 

\noindent {\bf CASE 2:} m=1.  Then $Q_{P_n \times \langle\sigma_n\rangle}$
is a contractible 1-dimensional complex with 3 vertices and two edges.
Define the maximal $p$-admissible tree $T$ of rank $n$ 
to be the tree with three
vertices $*$, $\circ$, $v$ and two edges $e_1$, $e_2$ where $e_1$ goes
from $*$ to $\circ$ and $e_2$ goes from $\circ$ to $v$.
The middle vertex of the 1-dimensional complex 
$Q_{P_n \times \langle\sigma_n\rangle}$
corresponds to the graph 
$\Gamma_T$. The other two vertices and two edges 
$Q_{P_n \times \langle\sigma_n\rangle}$ correspond 
to the two possible ways that $\Gamma_T$ can be collapsed
equivariantly.

\input{hyper2.pic}
\begin{figure}[here]
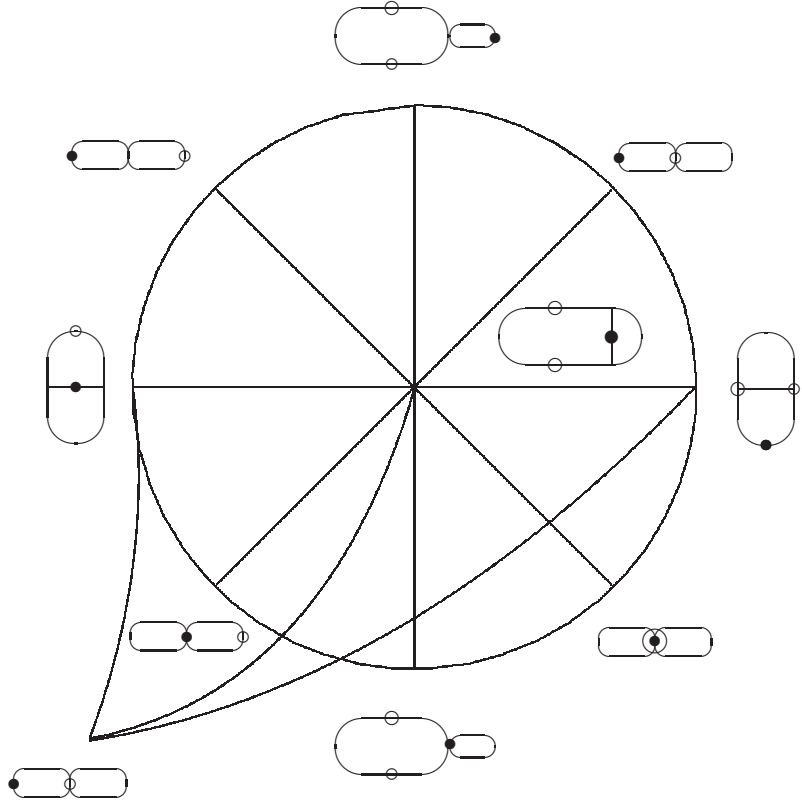

\caption{\label{hyper2} Simplices from the second maximal graph}
\end{figure}

\noindent {\bf CASE 3:} m=2. Then $Q_{P_n \times \langle\sigma_n\rangle}$
is a 2-dimensional complex with 13 vertices, 28 edges, and 16
two-simplices.  There are two maximal graphs in 
$Q_{P_n \times \langle\sigma_n\rangle}$.  Simplices coming from the
first graph are listed in figure \ref{hyper1} and simplices from 
the second graph are listed in figure \ref{hyper2}.  In figures
\ref{hyper1} and \ref{hyper2}, the maximal graphs are listed in the center.
These maximal graphs can be collapsed in various ways, and these are listed
around the periphery of the figures.  In the graphs, a solid dot
indicates the basepoint $*$ and the hollow dots represent attaching
points $\circ$ for the $\theta$-graph $\theta_{p-1}$.  If there is only one
hollow dot in a graph, both ends of the $\theta$-graph should be attached
to that one vertex.  Upon identifying the boundaries of the simplices 
listed in figures \ref{hyper1} and \ref{hyper2}, we obtain the
complex $Q_{P_n \times \langle\sigma_n\rangle}$ pictured in
figure \ref{hyper3}.  The complex is homeomorphic to the fletching
of a dart, three half disks, all identified along a
common line in their boundary.  
This complex is clearly contractible. \end{proof}

\centerline{\includegraphics{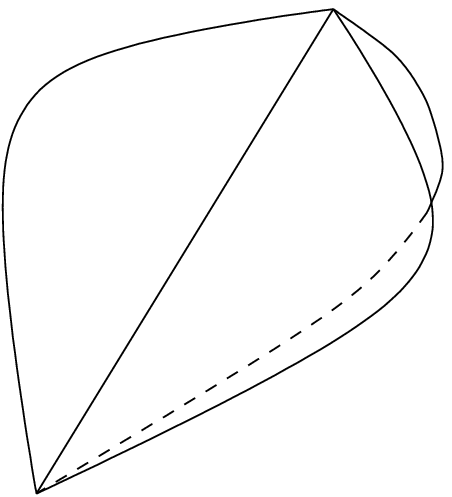}}
\begin{figure}[here]
\caption{\label{hyper3} The complete complex $Q_{P_n \times \langle\sigma_n\rangle}$}
\end{figure}

\end{document}